\documentclass{elsart}

\usepackage{amssymb, color, graphicx}
\usepackage[all]{xy}
\xyoption{dvips}

\newcommand{\abs}[1]{{\left\vert #1 \right\vert}}

\def\a{\alpha}
\def\b{\beta}
\def\g{\gamma}

\def\e{\epsilon}

\def\h{\eta}

\def\k{\kappa}

\def\x{\xi}

\def\s{\sigma}

\def\ph{\phi}

\def\w{\omega}
\def\G{\Gamma}
\def\D{\Delta}

\def\Om{\Omega}
\def\infinity{\infty}

\def\addots{\mathinner{\mkern1mu\raise1pt\vbox{\kern7pt\hbox{.}}\mkern2mu    \raise4pt\hbox{.}\mkern2mu\raise7pt\hbox{.}\mkern1mu}} 

\def\({\left(}
\def\){\right)}

\def\del{\partial}

\def\bo#1{{\bf #1}}

\def\RR{\bo R}
\def\ZZ{\bo Z}

\newcommand{\emb}{{\mathrm{Emb}}}
\newcommand{\imm}{{\mathrm{Imm}}}

\newcommand{\hofiber}{{\mathrm{hofiber}}}

\newcommand{\holim}{{\mathrm{holim}}}
\newcommand{\hocolim}{{\mathrm{hocolim}}}

\newcommand{\map}{{\mathrm{map}}}

\newcommand{\colim}{{\mathrm{colim}}}

\newcommand{\link}{{\mathrm{Link}}}
\newcommand{\conn}{{\mathrm{conn}}}

\begin{document}

\begin{frontmatter}



\title{A Manifold Calculus Approach to Link Maps and the Linking Number}


\author{Brian A. Munson}
\ead{\texttt{munson@math.stanford.edu}}
\address{Department of Mathematics, Stanford University, Stanford, CA 94305, USA}

\begin{abstract}
We study the space of link maps $\link(P_1,\ldots P_k; N)$, which is the space of maps $P_1\coprod\cdots\coprod P_k\rightarrow N$ such that the images of the $P_i$ are pairwise disjoint. We apply the manifold calculus of functors developed by Goodwillie and Weiss to study the difference between it and its linear and quadratic approximations. We identify an appropriate generalization of the linking number as the geometric object which measures the difference between the space of link maps and its linear approximation. Our analysis of the difference between link maps and its quadratic approximation connects with recent work of the author, and is used to show that the Borromean rings are linked.
\end{abstract}

\begin{keyword}
links \sep calculus of functors \sep linking number
\PACS
\end{keyword}
\end{frontmatter}

	
\section{Introduction}\label{intro}


Let $P_1,\ldots P_k$ be smooth compact, closed manifolds of dimensions $p_1,\ldots, p_k$, and let $N$ be a smooth manifold.

\begin{defn}
The space of \emph{link maps} of $P_1,\ldots, P_k$ in $N$ is the space $\link(P_1,\ldots,P_k;N)=\{f_i:P_i\rightarrow N | f_i(P_i)\cap f_j(P_j)=\emptyset\mbox{ for $i\neq j$}\}$, topologized as a subspace of the space of maps.
\end{defn}

It is important to distinguish this from the space of embeddings of the disjoint union, since all we assume about the individual maps $f_i$ is that they are smooth, so that we ignore self-linking. Alternately, one could study the case where the $f_i$ are immersions, a case considered by Hatcher-Quinn \cite{hq} and Schneiderman-Teichner \cite{st}. Our inspiration for studying the space of link maps was to get more insight into the space of embeddings. Our description of the linking number was inspired by section 1.4 of \cite{gkw2}, which describes double point obstructions to embeddings, and we also rely heavily on \cite{mun}, in which the author measured secondary obstructions (to the double point obstruction) for embedding a manifold in Euclidean space. 


The study of the space of link maps begins with the work of Milnor \cite{mil3}, who worked in the classical case where the $P_i$ are circles and $N$ is a $3$-manifold, and was mostly concerned with the fundamental group of the complement of a given link map. This was further investigated by Levine \cite{lev}, Habegger and Lin \cite{hablin}, and many others. The higher dimensional analog, which we are most concerned with here, has been extensively studied by Hatcher and Quinn \cite{hq}, Koschorke \cite{kos1,kos2,kos3,kos4}, Habegger and Kaiser \cite{hak}, Skopenkov \cite{sko1}, Schneiderman and Teichner \cite{st}, and more recently, Klein and Williams \cite{kw}.

Our analysis of ``linking manifolds'' in section \ref{linktot1} owes much to the work of Hatcher and Quinn \cite{hq}, as well as to the survey paper of Goodwillie, Klein, and Weiss \cite{gkw1}. Our work also intersects the very interesting recent work of Klein and Williams \cite{kw}. Section 9 of \cite{kw} specifically talks about linking and the linking number, and we would like to emphasize that our work on the linking number was done independently. Our work on the linking number also relates to the recent work of Chernov and Rudyak \cite{cr}. The difference is that, for us, a linking number depends on the choice of a path between the link in question and the unlink, while \cite{cr} wishes to keep track of such choices.

Our goal is to understand the space of link maps from the point of view of the manifold calculus of functors developed by Weiss and Goodwillie \cite{we1,gw}. We consider the space of link maps $\link(P_1,\ldots, P_k;N)$ as a contravariant functor from the poset of open subsets of $P_1\coprod \cdots\coprod P_k$. We describe briefly in Section \ref{manifoldcalculus} how this theory assigns $k^{th}$ degree ``polynomial'' approximations to such a functor $F$, denoted $\mathcal{T}_kF$. We will give a geometric description of the fiber of $\link(P_1,\ldots,P_k;N)\rightarrow\mathcal{T}_j\link(P_1,\ldots,P_k;N)$ when $j=1,2$. This geometric description is in terms of cobordism spaces, the relevant details of which can be found in \cite{mun}, and the relevant definitions and results of which are recounted in Section \ref{cobordismspaces}. Our Theorems \ref{t1mapofspaces} (resp. \ref{t2mapofspaces}) state that there is a map from the homotopy fiber of the space of link maps to its linear (resp. quadratic) approximation to a cobordism space. As the difference between the space of link maps and its linear (resp. quadratic) approximation begins with quadratic (resp. cubic) information, the relevant cobordism spaces are models for the quadratic (resp. cubic) homogeneous parts of the Taylor tower for these functors (see Section \ref{manifoldcalculus} for explanation of this terminology).

\begin{thm}\label{t1mapofspaces}
There is a map of spaces $$l_2:\hofiber(\link(P_1,P_2;N)\rightarrow\map(P_1,N)\times\map(P_2,N))\rightarrow\Om C_2(P_1,P_2;N).$$
\end{thm}

\begin{thm}\label{t2mapofspaces}
There is a map of spaces $$l_3:\hofiber(\link(P_1,P_2,P_3;N)\rightarrow\mathcal{T}_2\link(P_1,P_2,P_3;N))\rightarrow \Om C_{3}(P_1,P_2,P_3;N).$$
\end{thm}

One can ask whether or not these cobordism spaces describe these homotopy fibers in the sense that there is a highly connected map between them, and we conjecture that they do. Conjectures \ref{c2connectivity} and \ref{c3connectivity} state that the cobordism models we describe admit a highly connected maps from the difference between the space of link maps and the polynomial approximations given by the application of the manifold calculus, and we conjecture the range in which this connectivity holds. We should note that the connectivity estimates for spaces of embeddings, due to Goodwillie and Klein \cite{gk} are extremely difficult to obtain, and it is possible this will be the case for link maps as well. We will begin to address such connectivity statements in a separate paper \cite{mun2}. 

The two applications given for our Theorems \ref{t1mapofspaces} and \ref{t2mapofspaces} are to identify the map of spaces given in Theorem \ref{t1mapofspaces} with the linking number in the classical case, and to prove, using Theorem \ref{t2mapofspaces}, that the Borromean rings are linked.

Theorem \ref{t2mapofspaces} requires the use of manifolds with corners, and Section \ref{manifoldswithcorners} addresses the relevant background. Sections \ref{cobordismspaces} gives the necessary background for cobordism spaces. Section \ref{manifoldcalculus} gives the necessary background for manifold calculus and contains the models for $\mathcal{T}_1\link$ and $\mathcal{T}_2\link$ which we utilize. Section \ref{linktot1} is devoted to constructing the map of spaces appearing in Theorem \ref{t1mapofspaces}, and Section \ref{classicalcase} discusses why this is the same as the linking number in all of the classical cases. Section \ref{linktot2} is spent constructing the map of spaces in Theorem \ref{t2mapofspaces}, and we use this map in Section \ref{borromeanrings} to prove that the Borromean rings are linked.



\subsection{Conventions}\label{conventions}

We write $QX$ for $\Om^{\infty}\Sigma^{\infty}X$ where X is a based space, and $Q_+X$ for $Q(X_+)$ when $X$ is unbased, and where $+$ denotes a disjoint basepoint. When we say a map is an \textsl{equivalence}, we mean it is a weak equivalence, unless otherwise noted. For a vector bundle $\xi$ over a space X, we denote by $T(X;\xi)$ is Thom space. We will use the same capital $T$ for tangent bundles, but this should cause no confusion. If $X$ is a finite dimensional unbased space with vector bundles $\xi$ and $\eta$ over it, choose a vector bundle monomorphism $\eta\rightarrow \e^i$, and let $Q_+(X;\xi-\eta)$ be $\Om^iQ_+T(X;\xi\oplus\e^i/\eta)$. We write Spaces for the category of fibrant simplicial sets, and we work in this category unless otherwise noted. Thus $\map(M,N)$ is the simplicial set whose $k$-simplices are the fiber-preserving smooth maps of $M\times\Delta^k\rightarrow N\times\Delta^k$. By fiber-preserving we mean that if $f_k$ is a $k$-simplex of $\map(M,N)$ and $p_N:N\times\Delta^k\rightarrow\Delta^k$ is the projection, then the composition $p_N\circ f=p_M$, where $p_M:M\times\Delta^k\rightarrow\Delta^k$ is the projection. Other mapping spaces are translated to the category of simplicial sets in a similar manner.

\section{Manifolds with corners}\label{manifoldswithcorners}

A manifold with corners is a generalization of the notion of a manifold with boundary. The idea is that the boundary of a manifold with corners is allowed to have boundary, and this boundary is allowed to have boundary, and so on. The easiest example is the solid $n$-cube. It would be nice if we could extend the definition to include things like the boundary of this solid $n$-cube, but the problem is that in general there are many ways codimension $k$ submanifolds can meet when $k>1$. It is therefore cumbersome to make a good working definition, and this goes beyond our aim. We will nonetheless need to deal with the issue of giving a smooth structure to a manifold formed by gluing together several manifolds with corners in the formation of the obstruction manifold in Section \ref{linktot2}. We have two goals. The first is to explain how to make a smooth manifold by gluing together smooth manifolds with corners, at least in the case where the codimension of the corners is small. The second is to explain how to glue together maps of vector bundles, given over a collection of manifolds with corners which glue together to form a smooth manifold, to make a vector bundle map over their union. The former will be used to show that the manifold defined in Section \ref{linktot2} is smooth, and the latter to show that vector bundle isomorphisms, given by transversality, can be glued together to form a vector bundle isomorphism over that manifold. Our definition of a manifold with corners is a modification of the definition of a smooth manifold given by Milnor in \cite{mil2}. Following this work, we first discuss smooth functions on half-spaces.

Let $\RR^m_+=[0,\infty)^m$. We can also think of $\RR^m_+$ as a subset of $\RR^n$ for $n>m$ by the inclusion of $\RR^m$ in $\RR^n$ as the first $m$ coordinates. For $0\leq k\leq m$, let $\del_k\RR^m_{+}$ denote the subspace where at least $k$ of the coordinates are zero. We have $\del_m\RR^m_+\subset \cdots\subset \del_k\RR^m_+\subset \cdots\subset \del_0\RR^m_+=\RR^m_+$. We call $\del_k\RR^m_+$ the \emph{$k$-stratum} of $\RR^m_+$.

\begin{defn}
We say that a map $f:\RR^m_+\rightarrow \RR^n_+$ is \emph{smooth} if it is the restriction of some smooth function $f':\RR^m\rightarrow \RR^n$.
\end{defn}

\begin{defn}
A subset $M\subset \RR^k$ is a \textsl{smooth $m$-manifold with corners} if each $x\in M$ has a neighborhood $U\cap M$ that is diffeomorphic to an open subset $V\cap \RR^m_{+}$. The $k$-\emph{stratum} (or \emph{codimension $k$ boundary}) of $M$ is the set of all points in $M$ that correspond to points of $\del_k\RR^m_+$ under such a diffeomorphism, and we denote this set by $\del_k M$.
\end{defn}

Let $M$ is a smooth $m$-manifold with corners, and let $k\leq m$ be fixed. If $M$ requires only charts of the form $\RR^i_+\times\RR^{m-i}$ for $i\leq k$, then we say \emph{$M$ has at most $k$-strata}. The manifolds we will ultimately be interested in have at most $2$-strata.

\subsection{Tangent Space}

To define the tangent space of a smooth manifold $M^m\subset \RR^k$ at $x\in M$ we first pick a parametrization $g:U\rightarrow M$ of a neighborhood $g(U)$ of $x\in M$ with $g(u)=x$, and since $M\subset \RR^k$ we may think of this as a map $g:\RR^m\rightarrow \RR^k$. We then define the tangent space $T_xM$ to be the image of $dg_u$. Now if $M$ is a manifold with corners whose interior is a smooth manifold, we can define as above the tangent space to $M$ at any point. In this case the parametrization $g$ is a smooth map $g:\RR^m_+\rightarrow \RR^k$, and by definition $g$ is the restriction of a smooth map $\RR^m\rightarrow \RR^k$.

\subsection{Gluing and smooth structures}\label{gluemanifoldswithcorners}

We will only discuss gluing together manifolds with at most $2$-strata, as that is all we require. The case $k=1$ is simple and well understood. A manifold with corners with only $1$-strata is a manifold with boundary. Suppose one has two manifolds $M_1$ and $M_2$ with boundary and a diffeomorphism $j:\del M_1\rightarrow \del M_2$. Choose collars $\del M_1\times [0,1)\subset M_1$ and $\del M_2\times [0,1)\subset M_2$. We think of these as embeddings of $\del M_i\times [0,1)$ in $M_i$ such that $(x,0)$ maps to $x$. Then the union $M=M_1\cup_j M_2$ has the structure of a smooth manifold, since now every point in $M$ has a neighborhood homeomorphic to $\RR^m$. In particular, it is then clear what we mean by its tangent bundle $TM$. If we denote by $\del$ the subset of $M$ corresponding to $j(\del M_1)=\del M_2$, then what we have shown is that $\del$ has a neighborhood in $M$ homeomorphic to $\del \times B^1$, where $B^1$ is an open $1$-disk.

One can generalize the existence of collars to manifolds with corners in the following way.

\begin{prop}\label{stratumneighborhood}
Let $M^m$ be a smooth manifold with corners with $k$-stratum $\del_kM$. Then there exists a neighborhood $N$ of $\del_kM$ in $M$ and a diffeomorphism $N\cong\del_kM\times[0,1)^k$.
\end{prop}

One can mimic the proof using transversality Hirsch gives of the collaring theorem in \cite{hir}.  We are going to use this to give a smooth structure to a closed topological manifold which is made from the union of smooth manifolds with corners. For reasons we described above, we will only discuss manifolds with at most $2$-strata.

\begin{prop}\label{atmost2stratasmoothmanifold}
Let $M_1, \ldots, M_n$ be manifolds with corners with at most $2$-strata such that $\del_1M_i=\del_1^{i+1}M_i\cup\del_1^{i-1}M_i$ and $\del_1^{i+1}M_i\cap\del_1^{i-1}M_i=\del_2M_i$ for all $i$ (the subscripts are to be read modulo $n$). Let $j_{i,i+1}\del_{1}^{i+1}M_i\rightarrow\del_1^{i}M_{i+1}$ be diffeomorphisms, where the subscripting integers are read modulo $n$. Suppose that the composition $j_{n,1}\circ\cdots\circ j_{1,2}=id$ when restricted to $\del_2M_i$. Then the union $M=\cup_i M_i$ can be given the structure of a smooth closed manifold.
\end{prop}

The idea is to show that the image of the $2$-stratum $\del_2$ in $M$ should have a neighborhood homeomorphic to $\del_2\times B^2$, where $B^2$ is the open $2$-disk. 




\begin{pf}
Divide $\RR^2$ up into $n$ equal \emph{sectors} $S_i$ for $i=1, \ldots n$. By a sector we mean the area between two rays meeting at the origin in angle $2\pi/n$. Denote these rays by $\del_{i-1}S_i$ and $\del_{i+1}S_i$, for $i=1$ to $n$ modulo $n$.





Then $\del_2M_i$ has a neighborhood $U_i$ in $M_i$ diffeomorphic to $\del_2M_i\times [0,1)^2$ by proposition \ref{stratumneighborhood}. We then choose diffeomorphisms $\del_2M_i\times [0,1)^2\cong \del_2M_i\times S_i$ in compatible with the diffeomorphisms $j_{i,i+1}$. By compatible we mean that the following diagram should commute.

$$\xymatrix{
\del_1^{i+1}U_i\ar[dd]_-{\cong}\ar[rr]^-{j_{i,i+1}} & & \del_1^{i}U_{i+1}\ar[dd]^-{\cong}\\
 & & \\
\del_2\times\del_{i+1}S_i\ar[rr]_-{=} & & \del_2\times\del_iS_{i+1}\\
}$$


Then it follows that there is a neighborhood of $\del_2$ in $M$ homeomorphic to $\del_2\times B^2$, and we use this homeomorphism to give $M$ its smooth structure. This gives $M$ a tangent bundle. Moreover, at each point in $\del_2$, the restriction of the tangent bundle to each sector $S_i$ is the tangent bundle already given to $\del_2$ on that sector by the inclusion of $M_i$.
$\Box$\end{pf}

\subsection{Vector bundles}

\begin{prop}\label{gluebundles1strata}
Suppose we have manifolds $M_1$ and $M_2$, continuous proper maps $f_i:M_i\rightarrow X$, and isomorphisms $\phi_i:TM_i\oplus f_i^*(\xi)\rightarrow f_i^*(\eta)$ for $i=1,2$, and a diffeomorphism $j:\del M_1\rightarrow \del M_2$ satisfying $f_1 = f_2\circ j$ on the boundary. Suppose that the restriction of $\phi_1$ and $j^\ast\phi_2$ to $\del M_1$ and $\del M_2$ respectively are homotopic. Then $M=M_1\cup_j M_2$ is a smooth manifold and there is an isomorphism $\phi:TM\oplus f^*(\xi)\rightarrow f^*(\eta)$ compatible with $\phi_i$ for $i=1,2$.
\end{prop}

\begin{pf}
Using collars of $\del M_i$ and the diffeomorphism $j$ we can make a smooth manifold $M=M_1\cup_j M_2$ with a continuous map $f:M\rightarrow X$, as discussed in section \ref{gluemanifoldswithcorners}. The only thing left is to make a commutative diagram


$$\xymatrix{
TM_1|_{\del M_1}\oplus f_1^*(\xi)|_{\del M_1} \ar[rr]^-{\phi_1|_{\del M_1}}\ar[dd] & & f_1^*(\eta)|_{\del M_1}\ar[dd] \\
 & & \\
j^*TM_2|_{\del M_2}\oplus j^*f_2^*(\xi)|_{\del M_2} \ar[rr]^-{j^*\phi_2|_{\del M_2}} & & j^*f_2^*(\eta)|_{\del M_2}\\
}
$$

where the horizontal isomorphisms are the $\phi_i$, the right vertical isomorphism is the identity, and the leftmost isomorphism is given as follows. We write $TM_i|_{\del M_i}=T\del M_i\oplus \e$ by identifying $\e$ with $\nu(\del M_i\subset M_i)$ and letting $1$ correspond to the outer unit normal via some Riemannian metric on $M_1$, and let $1$ correspond to the inner unit normal on $M_2$. Then the isomorphism between $T\del M_1\oplus \e$ and $j^*T\del M_2\oplus j^*\e$ is the obvious one (it is $1$ on the $\e$ part). If this diagram commutes, then the isomorphisms $\phi_i$ agree on $\del M_i$, so we have produced a definite isomorphism $\phi:TM\oplus f^*(\xi)\rightarrow f^*(\eta)$ at all points $x\in M$. In fact, it is enough if the diagram above commutes up to homotopy. Let $\phi_t$ be the homotopy from $\phi_1$ to $\phi_2$. For notational convenience, let us suppress the diffeomorphism $j$, and denote $\del M_i$ by simply $\del$. We have $M=\colim(M_1\leftarrow \del \rightarrow M_2),$ and we have bundle isomorphisms on $M_i$ and a homotopy between them on $\del$. If we set $M'=\hocolim(M_1\leftarrow\del\rightarrow M_2),$ then we have over each point in $M'$ a definite isomorphism $TM\oplus f^*(\xi)\cong f^*(\h)$. It is $\phi_1$ on $M_1$, $\phi_2$ on $M_2$, and $\phi_t$ on $\del\times \{t\}$. The canonical map $M'\rightarrow M$ is a homotopy equivalence, and we can pull back the bundle isomorphism on $M'$ to $M$ by a homotopy inverse to get the desired bundle isomorphism.
$\Box$\end{pf}

\begin{rem}
In our case the bundle isomorphisms are going to be given by transversality. That is, they will be induced by the derivatives of certain functions which define manifolds with corners, and we will need to check that they are homotopic. 

As the proof of Proposition \ref{gluebundles1strata} indicates, the issues of giving a smooth structure to the union $M_1\cup M_2$ and that of gluing together the bundle isomorphisms can be dealt with separately. We next consider the process of forming a map of vector bundles over a space which is the union of other spaces by gluing together vector bundle maps given on the smaller pieces.
\end{rem}


Let $\a$ and $\b$ be vector bundles over a space $Z$. 

\begin{defn}\label{maxrankbundlemaps}
Define $L_{\max}(\a,\b)_Z$ to be the space of vector bundle maps from $\a$ to $\b$ over $Z$ such that the linear map of fibers is of maximal rank.
\end{defn}

We generally suppress $Z$ from this notation and just write $L_{\max}(\a,\b)$, because the space in question will be clear from the context. Our immediate goal is to describe a generalization of Proposition \ref{gluebundles1strata}. This is contained in the next two propositions.

Let $\mathcal{C}$ be a small category, and suppose we are given a covariant functor $\mathcal{X}:\mathcal{C}\rightarrow\mbox{ Spaces}$, whose values we will denote $X(c)$, and let us assume that $X(c)$ is compact. Let $X(\mathcal{C})$ denote the diagram of spaces obtained from $\mathcal{X}$. Suppose in addition that $Y$ is some fixed space, and $\xi$ and $\eta$ are vector bundles over $Y$. Let $f_c:X(c)\rightarrow Y$ be maps for each $c\in\mathcal{C}$ such that the diagram $X(\mathcal{C})\rightarrow Y$ commutes. Then $X=\hocolim_{\mathcal{C}}\mathcal{X}$ is a space with a map $f:X\rightarrow Y$ given by $f_c$.

The model for the homotopy colimit of a functor $F:\mathcal{C}\rightarrow \mbox{ Spaces}$ we have in mind is the coequalizer

$$
\xymatrix{
\hocolim_\mathcal{C}F={\mbox{coeq}}(\coprod_{c\rightarrow c'}\abs{c'\downarrow\mathcal{C}}\times F(c)\ar@{->}[r]<3pt>\ar@{->}[r]<-3pt> &\coprod_{c}\abs{c\downarrow\mathcal{C}}\times F(c))
}
$$

\begin{prop}\label{gluebundles1}
Let $\mathcal{C}$, $\mathcal{X}$, $Y$, $\xi$, $\eta$, and $\{f_c\}_{c\in\mathcal{C}}$ be as above. Consider the functor $\mathcal{E}:\mathcal{C}\rightarrow\mbox{ Spaces}$ for which $\mathcal{E}(c)=L_{\max}(f_c^\ast\xi,f_c^\ast\eta)$ as a bundle over $X(c)$ and the maps $\mathcal{E}(c)\rightarrow\mathcal{E}(c')$ are given by pullback of the map $X(c)\rightarrow X(c')$. Then $E=\hocolim_{\mathcal{C}}\mathcal{E}$ is a vector bundle over $X=\hocolim_{\mathcal{C}}\mathcal{X}$.
\end{prop}

\begin{pf}
This follows from the fact that the maps $\mathcal{E}(c)\rightarrow\mathcal{E}(c')$ are given by the pullback by $X(c)\rightarrow X(c')$. $\Box$
\end{pf}

We wish to describe how to build a section of $E\rightarrow X$ from sections of $\mathcal{E}(c)\rightarrow\mathcal{X}(c)$. Let $\chi_{c\rightarrow c'}:X(c)\rightarrow X(c')$ and $\Xi_{c\rightarrow c'}:L_{\max}(f_c^\ast\xi,f_c^\ast\eta)\rightarrow L_{\max}(f_{c'}^\ast\xi,f_{c'}^\ast\eta)$ denote the maps given by the functors $\mathcal{X}$ and $\mathcal{E}$ respectively. Let $u_{c\rightarrow c'}:\abs{c'\downarrow \mathcal{C}}\rightarrow\abs{c\downarrow \mathcal{C}}$. We will suppress the morphism $c\rightarrow c'$ from these notations unless it would cause confusion.

\begin{prop}\label{gluebundles2}
Suppose we are given maps $\Phi_{c\rightarrow c'}:\abs{c'\downarrow \mathcal{C}}\times X(c)\rightarrow L_{\max}(f_c^\ast\xi,f_c^\ast\eta)$. If these maps satisfy $\Phi(s',x)=\Phi(u(s'),x)$ and $\Xi(\Phi(s,x)=\phi(s,\chi(x))$, then they piece together to form a section of $E\rightarrow X$.
\end{prop}

\begin{pf}
This follows immediately from the coequalizer definitions of $E$ and $X$.$\Box$
\end{pf}

We will be interested in the case when then $X(c)$ are manifolds with corners, $X:=\colim_{\mathcal{C}}X$ has the structure of a smooth manifold, and the map $X'\rightarrow X$ is a homotopy equivalence so that we can pull back the bundle map obtained over $X'$ by a homotopy inverse to $X$. In our case, the maps $\Phi_{c\rightarrow c'}$ will be given by various homotopies, much as in Proposition \ref{gluebundles1strata}. See Lemma \ref{piecetogetherbundledata}.


\section{Cobordism spaces}\label{cobordismspaces}

We begin with a very brief description of cobordism spaces which we will employ in our description of $\mathcal{T}_i\link(P_1,\ldots,P_k;N)$ for $i=1,2$. These spaces were used extensively in \cite{mun}. Identifying a good model for these spaces is necessary to define the maps in Theorems \ref{linktot1} and \ref{linktot2}.\\


Let $X$ be a space, and $\x$ and $\h$ vector bundles on $X$. An element of the cobordism group $\Om_{k}^{\x-\h}(X)$ is represented by a triple $(W^k,f,\ph)$ (sometimes denoted by just $W$) where $W$ is a $k$-dimensional smooth manifold embedded in $\RR^{\infinity}$, $f:W \rightarrow X$ is continuous and proper, and $\ph$ is a stable isomorphism $TW \oplus f^*\x \cong f^*\h$. The equivalence relation for representatives is the usual one defined by $(k+1)$-dimensional manifolds with boundary. We seek a space whose homotopy groups are the cobordism groups described above, and we call such a space a cobordism space.

\begin{defn}[Simplicial Model for a Cobordism Space]\label{cobspacedefn} The simplicial set $C_{\bullet}^{\xi-\h}(X)$ has as its $k$-simplices the set $C_k=\{(W^{d+k},f,\ph)\}$ where $d=\dim(\h)-\dim(\x)$, $W$ is a smooth $(k+d)$-dimensional manifold embedded in $\RR^{\infty}\times \D^k$, $W$ is transverse to $\RR^{\infty}\times \del_S\D^k$ for all nonempty subsets $S\subset\{0,1,\ldots, k\}$, $f:W\rightarrow X$ is continuous and proper, and $\ph: TW\oplus f^*(\xi)\rightarrow f^*(\eta)$ is a stable isomorphism.
\end{defn}

The manifolds $W^{d+k}\subset \Delta^k\times\RR^\infty$ should the \emph{conditioned}. To be conditioned means that if we denote by $W_t$ the part of $W$ that sits over $t\in\Delta^k$, then $W_t$ should be independent of $t$ in a neighborhood of $\cup_i \del_i\Delta^k$.\\


The face and degeneracy maps are induced by those of $\D^{\bullet}$. The $i^{\mbox{th}}$ face map $d_i:C_k\rightarrow C_{k-1}$ is just the intersection of $W^{d+k}$ with the $i^{\mbox{th}}$ face of $\D^k$. The $i^{\mbox{th}}$ degeneracy map $s_i:C_k \rightarrow C_{k+1}$ takes $W$ to the fiber product $W'$

$$\xymatrix{
W' \ar[r]\ar[d] & \RR^{\infty}\times \D^{k+1}\ar[d]^{s_i} \\
W\ar[r] & \RR^{\infty}\times \D^{k}\\
}
$$

where $s_i$ is the $i^{th}$ degeneracy for $\D^{\bullet}$. That $C_{\bullet}^{\xi-\h}(X)$ satisfies the axioms for a simplicial set is straightforward, because we are building on the usual simplicial structure on $\Delta^k$.\\

The next three propositions are established in \cite{mun}.

\begin{prop}\label{CKan}
$C_{\bullet}^{\xi-\h}(X)$ is a Kan complex.
\end{prop}

\begin{prop}\label{loopcobordism}
There is an equivalence $C_{d+l}^{\xi-\h}(X)\simeq \Om^lC_{d}^{\xi-\h}(X).$
\end{prop}

\begin{rem}\label{cobspacethomspace}
This cobordism space is equivalent to $QT(X;\x-\h)$. To see the equivalence, consider the subcomplex of the total singular complex of $QT(X;\x-\h)$ consisting of those $k$-simplices $\k:\Delta^k\rightarrow \Om^n\Sigma^n(T(\x-\h))$ that correspond to maps $\k':\Sigma^n(\Delta^k)\rightarrow \Sigma^n(T(\x-\h))$ which are transverse to the zero section of $T(\x-\h)$. This sub-complex is equivalent to the full complex and the map $\k\mapsto \k'^{-1}(0)$ to the cobordism model is an equivalence. See \cite{tg1} for a similar construction. 
\end{rem}

That $C_\bullet$ is a Kan complex ensures that the homotopy groups of its realization will be the cobordism groups we want. That is, $\pi_k\abs{C_{d}^{\xi-\h}(X)}=\Om_{d+k}^{\xi-\h}(X)$. The second proposition, together with the next, will be useful in explicitly identifying these cobordism groups in special cases, and will be used in our identification of the linking number in Section \ref{classicalcase} and proof that the Borromean rings are linked in Section \ref{borromeanrings}.

%

\begin{prop}\label{countclasses}
If $X$ is connected, the group $\Om_0^{\xi-\eta}(X)$ is isomorphic with $\ZZ$ if $w(\xi)=w(\eta)$, and $\ZZ/2$ if $w(\xi)\neq w(\eta)$.
\end{prop}


\section{Manifold Calculus}\label{manifoldcalculus}

Manifold calculus, developed by Goodwillie and Weiss \cite{gw,we1} studies contravariant functors $F:\mathcal{O}(M)\rightarrow\mbox{ Spaces}$, where $\mathcal{O}(M)$ is the poset of open subsets of a smooth manifold $M$. Examples include $U\mapsto\emb(U,N)$, the space of embeddings of $U$ in a smooth manifold $N$, $U\mapsto\map(U,X)$, the space of maps of $U$ to a space $X$, and $U\mapsto\imm(U,N)$, the space of immersions of $U$ in a smooth manifold $N$. To such a functor $F$, the manifold calculus associates a tower of functors

$$
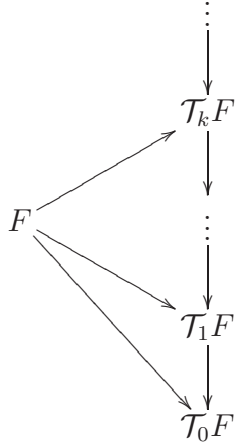
\xymatrix{
&& \vdots\ar[d]\\
&& \mathcal{T}_{k}F\ar[d]\\
F\ar[rru]\ar[rrd]\ar[rrdd] && \vdots\ar[d]\\
&& \mathcal{T}_1F\ar[d]\\
&& \mathcal{T}_0F\\
}
$$

called the \emph{Taylor tower} of $F$.

\begin{defn}
For a contravariant functor $F:\mathcal{O}(M)\rightarrow\mbox{ Spaces}$ we define the \textsl{$k^{th}$ Taylor approximation} to $F$, denoted $\mathcal{T}_kF:\mathcal{O}(M)\rightarrow\mbox{ Spaces}$, by

$$\mathcal{T}_kF(U)=\holim_{V\in \mathcal{O}_k(U)}F(V).$$

Here $\mathcal{O}_k(U)$ is the subcategory of $\mathcal{O}(U)$ consisting of those open sets $V\subset U$ which are diffeomorphic to at most $k$ open balls.
\end{defn}

\begin{defn}
We say that $F$ is \emph{polynomial of degree $\leq k$} if given pairwise disjoint closed subsets $A_0,A_1,\ldots,A_k$ of $U\in\mathcal{O}(M)$, the map $F(U)\rightarrow\holim_{\emptyset\neq S\subset \{0,1,\ldots, k\}}$ is a weak equivalence. We say $F$ is \emph{homogeneous of degree $k$} if additionally $\mathcal{T}_{k-1}F\simeq\ast$.
\end{defn}

For instance, the functors $U\mapsto\map(U,X)$ and $U\mapsto\imm(U,N)$ are polynomial of degree $\leq 1$, and the functor $U\mapsto\map(U^k,X)$ is a polynomial of degree $\leq k$. See \cite{we1} for proofs. The next two theorems state that the functors $\mathcal{T}_kF$ are polynomial and that they are essentially determined by their values on special open sets. 

\begin{thm}[\cite{we1}, Theorem 6.1]\label{tkpoly}
The cofunctor $\mathcal{T}_kF$ is polynomial of degree $\leq k$.
\end{thm}

\begin{rem}\label{sectionslinear}
For a smooth manifold $X$ and a vector bundle $\xi\rightarrow X$, let $\Gamma(\xi, X)$ denote the space of sections. $\Gamma(\xi,-)$ is a polynomial of degree $\leq 1$ from $\mathcal{O}(X)\rightarrow\mbox{ Spaces}$. Hence, $\Gamma(\xi,X)\simeq\holim_{U\in\mathcal{O}(X)}\Gamma(\xi,U)$. In order to produce a section defined on all of $X$ is is therefore enough to produce an open cover $\mathcal{U}$ of $X$, an element of $\Gamma(\xi,U)$ for each $U\in\mathcal{U}$, and homotopies between these sections on their intersections.
\end{rem}

\begin{thm}[\cite{we1}, Theorem 5.1]\label{polydet}
Suppose that $\g:F_1\rightarrow F_2$ is a morphism of good cofunctors, and that $F_i$ is polynomial of degree $k$ for $i=1,2$. If $\g:F_1(V)\rightarrow F_2(V)$ is a homotopy equivalence for all $V\in \mathcal{O}_k(M)$, then it is a homotopy equivalence for all $V\in \mathcal{O}(M)$.
\end{thm}

From its definition we see that the values of $\mathcal{T}_kF$ are completely determined by its values on $\mathcal{O}_k(M)$, so Theorem \ref{polydet} is not too surprising.

Finally, we state the classification theorem for homogeneous functors.

\begin{thm}[\cite{we1}, Theorem 8.1]\label{homogclass}
Let $F$ be a homogeneous cofunctor of degree $k$. Then there is an equivalence

$$F(V)\rightarrow \G^c\left(p;{V\choose k}\right)$$

where $V\in\mathcal{O}(M)$, and $\G^c\left(p;{V\choose k}\right)$ is the space of compactly supported sections of a fibration $p:E\rightarrow {V\choose k}$.
\end{thm}

\subsection{The categories $\mathcal{O}_l(P_1\coprod \cdots\coprod P_k)$}

Let $P_1, \ldots P_k$ be smooth closed compact manifolds of dimension $p_1,\ldots p_k$. We wish to apply manifold calculus to study the space $\link(P_1,\ldots P_k; N)$. The functor in question is $\link(-;N):\mathcal{O}(P_1\coprod\cdots\coprod P_k)\rightarrow\mbox{ Spaces}$.

\begin{prop}\label{opensetsofdisjointunion}
There is an equivalence of categories $\mathcal{O}_l(P_1\coprod \cdots\coprod P_k)\simeq \coprod_{\sum_{j=1}^li_j =k} \mathcal{O}_{i_1}(P_1)\times\cdots\times\mathcal{O}_{i_k}(P_k)$ which sends an open set $U$ to $(U\cap P_1,\ldots, U\cap P_k)$.
\end{prop}

\begin{prop}\label{holimoverdisjointunion}
Let $\mathcal{E}$ be a category which is the disjoint union of two categories $\mathcal{C}$ and $\mathcal{D}$, whose morphism set is the disjoint union of the morphism sets of $\mathcal{C}$ and $\mathcal{D}$. Then for a functor $F:\mathcal{E}\rightarrow Top$, $\holim_{\mathcal{E}}F=\holim_{\mathcal{C}}F\times\holim_{\mathcal{D}}F$.
\end{prop}

We end this section by giving mapping space models for $\mathcal{T}_1\link(P_1,\ldots, P_k;N)$ and $\mathcal{T}_2\link(P_1,\ldots, P_k;N)$ that we will use.

\begin{prop}\label{t1model}
There is an equivalence $\mathcal{T}_1\link(P_1,\ldots ,P_k; N)\simeq \prod_{i=1}^k\map(P_i,N))$.
\end{prop}

\begin{pf}
Let $P=\coprod_i P_i$, $U\in\mathcal{O}(P)$, and set $U_i=U\cap P_i$. By definition and Propositions \ref{opensetsofdisjointunion} and \ref{holimoverdisjointunion}, 

\begin{eqnarray*}
\mathcal{T}_1\link(U_1,\ldots U_k;N)&\simeq&\prod_{i=1}^k\holim_{V_i\in\mathcal{O}_1(U_i)}\link(\emptyset, \ldots, V_i, \ldots, \emptyset;N)\\
 &\simeq&\prod_{i=1}^k\map(U_i,N).
\end{eqnarray*}
$\Box$\end{pf}

Alternatively, we could have used Theorem \ref{polydet} to prove this.

%

\begin{prop}\label{t2model}
The functor $\mathcal{T}_2\link(P_1,P_2;N)$ is equivalent to the homotopy pullback of 
$$
\xymatrix{
&\map(P_1\times P_2,N\times N-\Delta_{N})\ar[d]\\
\map(P_1,N)\times\map(P_2,N)\ar[r]&\map(P_1\times P_2,N\times N)\\
}
$$
\end{prop}



\begin{pf}
Both are polynomials of degree less than or equal to $2$, and the map from $\mathcal{T}_2\link(P_1,P_2;N)$ to the homotopy limit is an equivalence by inspection when $P_1\coprod P_2$ is a manifold consisting of at most two points. The result follows from Theorem \ref{polydet}. Compare Theorem 1.3 of \cite{gkw1} in the case $k=2$.
$\Box$\end{pf}

%

\begin{cor}
The functor $\mathcal{T}_2\link(P_1,\ldots, P_k;N)$ is equivalent to the homotopy pullback of 

$$\prod_{i<j}
\xymatrix{
&\map(P_i\times P_j,N\times N-\Delta_{N})\ar[d]\\
\map(P_i,N)\times\map(P_j,N)\ar[r]&\map(P_i\times P_j,N\times N)\\
}
$$
\end{cor}

\begin{pf}
This follows from Proposition \ref{holimoverdisjointunion}
$\Box$\end{pf}

Hence, we may think of an element of $\mathcal{T}_2\link(P_1,\ldots, P_k;N)$ as a tuple $(\bf{f},\bf{H}, \bf{F})$, where $\bf{f}=(f_1,\ldots,f_k)$ satisfies $f_i\in\map(P_i,N)$, $\bf{F}=(F_{12}, \ldots, F_{(k-1)k})$ satisfies $F_{ij}\in\map(P_i\times P_j,N\times N-\Delta_{N})$, and $H_{ij}$ is a homotopy between $F_{ij}$ and $f_i\times f_j$.

\section{Linking Manifolds and Quadratic Obstructions}\label{linktot1}

In this section we are going to examine the difference between $\link(P_1,P_2;N)$ and $\mathcal{T}_1\link(P_1,P_2;N)$. It should be clear from what follows that for $\link(P_1,\ldots, P_k;N)$ we make $k\choose 2$ linking manifolds, one for each pair $(i,j)$. Recall from Proposition \ref{t1model} that $\mathcal{T}_1\link(P_1,P_2;N)\simeq\map(P_1,N)\times\map(P_2,N)$.


Let $(f_1,f_2)\in\map(P_1,N)\times\map(P_2,N)$ be the basepoint, and let $$\a\in\hofiber(\link(P_1,P_2;N)\rightarrow\map(P_1,N)\times\map(P_2,N)).$$ This is a map $\a:I\rightarrow \map(P_1,N)\times\map(P_2,N)$ such that $\a(0)=(f_1,f_2)$, and $\a(1)\in\link(P_1,P_2;N)$. We will write $\a(t)=(f_{1,t},f_{2,t})$ so that $f_{1,0}=f_1$ and $f_{2,0}=f_2$.\\

\subsection{Cobordism space model for $\hofiber(\link(P_1,P_2;N)\rightarrow\mathcal{T}_1\link(P_1,P_2;N))$}

For the basepoint $(f_1,f_2)\in\map(P_1,N)\times\map(P_2,N)$, the images of $f_1$ and $f_2$ need not be disjoint for what we are about to do, and in general they might not be. We are mostly interested in the case when $(f,g)$ is in the image of $\link(P_1,P_2;N)$, but for what follows, this is not necessary. The constructions in this section were inspired by a very similar construction in section 1.4 of \cite{gkw1}. Compare also \cite{kw}.

\begin{defn}
$E_{12}=\holim(P_1\stackrel{f_1}{\rightarrow} N\stackrel{f_2}{\leftarrow} P_2)$
\end{defn}

Hence $E_{12}=\{(x_1,x_2,\w) : (x_1,x_2)\in P_1\times P_2, \w:[-1,1]\rightarrow N, \w(-1)=f_1(x_1),\w(1)=f_2(x_2)\}$

$E_{12}$ has maps to $P_1\times P_2$ and $N$ given by projection and evaluation of $\w$ at $0$. Hence we may pullback $TP_1\times TP_2$ and $TN$ to $E_{12}$. From this data we can make a cobordism space $C_{\bullet}^{TN-TP_1\times TP_2}(E_{12})$. For us, $d=p_1+p_2-n$, so that a $0$-simplex in this space is a manifold of dimension $(p_1+p_2-n)$. We will abbreviate this space by $C_2(P_1,P_2;N)$. It is equivalent to $QT(E_{12};TN-TP_1\times TP_2)$ by the Pontryagin-Thom construction.

We may assume that $f_1\times f_2:P_1\times P_2\rightarrow N\times N$ is transverse to $\Delta_N$, because the subcomplex of $\map(P_1,N)\times\map(P_2,N)$ of maps $(f_1,f_2)$ such that $f_1\times f_2$ is transverse to $\Delta_N$ is homotopy equivalent to the full complex.


\begin{defn}
Let $D=(f_1\times f_2)^{-1}(\Delta_N)$
\end{defn}

\begin{lem}
If $f_1\times f_2$ is transverse to $\Delta_N$, then $D$ is a smooth compact closed manifold of dimension $(p_1+p_2-n)$ with normal bundle $TN-TP_1\times TP_2$.
\end{lem}

\begin{pf}
This follows from transversality.
$\Box$\end{pf}

There is an inclusion $D\rightarrow E_{12}$ which associates the constant path $\w$ with each pair $(x_1,x_2)\in D$. Hence $D$ determines a point in $C_2(P_1,P_2;N)$.

We are now going to produce a path in $C_2(P_1,P_2;N)$ from our choice of $\a\in\hofiber(\link(P_1,P_2;N)\rightarrow\map(P_1,N)\times\map(P_2,N))$. Recall from above the notation $\a_t=(f_{1,t},f_{2,t})$. We may regard $f_{1,t}\times f_{2,t}:P_1\times P_2\rightarrow N\times N$ as a map $F:P_1\times P_2\times I\rightarrow N\times N$. We may assume that $F$ is transverse to the diagonal $\Delta_N\subset N\times N$, and that its restriction to $P_1\times P_2\times\{0\}$ and $P_1\times P_2\times\{1\}$ are both transverse to the diagonal. Thus we produce a smooth closed manifold $D_\a=F^{-1}(\Delta_N)\subset P_1\times P_2\times I$. Its normal bundle is isomorphic with $TN-TP_1\times TP_2$ by transversality, and $D_\a\subset E_{12}\times [0,1]$ is the subset of all $(x_1,x_2,\w,t)$ such that $f_{1,t}(x_1)=f_{2,t}(x_2)$, and $\w$ is the path 

$$\w(s)= \left\{ \begin{array}{ll}
                        f_{1,(1+s)t} & \mbox{if $-1\leq s \leq 0$}\\
                        f_{2,(1-s)t} & \mbox{if $0\leq s \leq 1$}
                    \end{array}
\right.$$

$D_\a$ gives us a nullcobordism of $D$, which is a cobordism between the empty manifold to $D$ in $C_2(P_1,P_2;N)$. Let $\Om_DC_2(P_1,P_2;N)$ denote the space whose $0$-simplices are nullcobordisms of $D$. To justify this notation note that if the images of $f$ and $g$ were disjoint, $D$ would be empty and the relevant space would be the loopspace $\Om C_2(P_1,P_2;N)$. This works for arbitrary maps of a $k$-simplex into $\hofiber(\link(P_1,P_2;N)\rightarrow\map(P_1,N)\times\map(P_2,N))$, producing a $k$-simplex in $\Om_DC_2(P_1,P_2;N)$, or equivalently a $k$-simplex in the space of paths from the basepoint of $QT(E_{12};TN-TP_1\times TP_2)$ to the image of $D$ by the Pontryagin-Thom construction. We need to note that the subcomplex of maps which satisfy the necessary transversality condition to define $D$ are homotopy equivalent to the full complex (see, for instance Hypothesis 3.18 of \cite{tg1}). This proves Theorem \ref{t1mapofspaces}, restated in a slightly different form here for convenience.

\begin{thm}
There is a map of spaces $$l_D:\hofiber(\link(P_1,P_2;N)\rightarrow\map(P_1,N)\times\map(P_2,N))\rightarrow\Om_DC_2(P_1,P_2;N).$$
\end{thm}

When the basepoint is in the image of the space of link maps, $D$ is empty and we write $l_2$ in place of $l_D$. The reason this map is important is that it is has a certain connectivity, depending on the dimensions of the manifolds involved.

\begin{conj}\label{c2connectivity}
$l_2$ is $(2n-\max\{2p+q,2q+p\}-3)$-connected.
\end{conj}

This follows from the following conjecture.

\begin{conj}\label{linktot2connectivity}
The map $\link(P_1,P_2;N)\rightarrow\mathcal{T}_2\link(P_1,P_2;N)$ is $(2n-\max\{2p+q,2q+p\}-3)$-connected.
\end{conj}

Here is a proof of Conjecture \ref{c2connectivity} using Conjecture \ref{linktot2connectivity} and results from Section 9 of \cite{kw}.

\begin{pf}
Let $F_l=\hofiber(\link(P_1,P_2;N)\rightarrow\mathcal{T}_1\link(P_1,P_2;N))$, and let $F_t=\hofiber(\mathcal{T}_2\link(P_1,P_2;N)\rightarrow\mathcal{T}_1\link(P_1,P_2;N))$. Consider the sequence of spaces $F_l\rightarrow F_t\rightarrow \Om C_3(P_1,P_2;N)$. Theorem B of \cite{kw} shows that the latter map is $(2n-p-q-3)$-connected, and Conjecture \ref{linktot2connectivity} implies the composed map has the desired connectivity.
$\Box$\end{pf}

\begin{rem}
The main difficulty in proving Conjecture \ref{linktot2connectivity} lies in the fact that the $P_i$ are not necessarily embedded in $N$. The difficulty is best expressed by the fact that the restriction map $\link(K_2;N)\rightarrow\link(K_1;N)$ for compact sets $K_1\subset K_2\subset P_1\coprod P_2$ is not in general a fibration. This makes it difficult to analyze homotopy fibers of such restrictions, unlike the case for embeddings. The simplest possible example is if we let $P_1=\{p\}$ is a single point, $P_2=\{q_1,q_2\}$, $K_1=P_1\coprod\{q_1\}$, and $K_2=P_1\coprod P_2$. The author will tackle Conjecture \ref{linktot2connectivity} and other similar statements in \cite{mun2}.
\end{rem}

\subsection{Properties of the Linking Manifold}

Now we consider the case where the basepoint $(f_1,f_2)$ is in the image of $\link(P_!,P_2;N)$ in $\map(P_1,N)\times\map(P_2,N)$. In this case we get a map 

\begin{equation}
\hofiber(\link(P_1,P_2;N)\rightarrow\map(P_1,N)\times\map(P_2,N))\rightarrow\Om C_2(P_1,P_2;N).
\end{equation}

As we have discussed above, for each $\a\in\hofiber(\link(P_1,P_2;N)\rightarrow\map(P_1,N)\times\map(P_2,N))$, we write $\a_t=(f_{1,t},f_{2,t})$ and use this to produce an element of $\Om C_2(P_1,P_2;N)$, which is represented by a smooth closed compact manifold $L$ of dimension $(p_1+p_2-n+1)$, called the \emph{linking manifold} of $\a$, and it plays the role of the linking number of $(f_1,f_2)$.




\begin{rem}
Clearly the linking manifold depends both on the choice of basepoint $(f_1,f_2)$ and the choice of $\a\in\hofiber(\link(P_1,P_2;N)\rightarrow\map(P_1,N)\times\map(P_2,N))$. The dependence on the basepoint is unavoidable, because we need to declare in advance what the unlink is, but in many cases, the linking manifold does not depend on the choice of path $\a$, only its values at $0$ and $1$. See section \ref{dependsonlift} for more.
\end{rem}

\begin{prop}
$L$ is a smooth closed compact $(p_1+p_2+1-n)$-dimensional manifold with stable normal bundle $TN-TP_1\times TP_2$, and defines an element of $\Om QT(E_{12};TN-TP_1\times TP_2)$.
\end{prop}

\begin{pf}
This follows from the existence of the map in Theorem \ref{t1mapofspaces}
$\Box$\end{pf}







\subsection{Dependence on lift}\label{dependsonlift}

In this section we discuss the extent to which the linking manifold depends on the choice of path $\a$ from $(f_1,f_2)$ to $(f_{1,1},f_{2,1})$. These results are simple but useful corollaries of Theorem \ref{t1mapofspaces}.

\begin{prop}\label{homotopiccobordant}
Suppose $\a,\b\in\hofiber(\link(P_1,P_2;N)\rightarrow\map(P_1,N)\times\map(P_2,N))$ and that $\a(1)=\b(1)$. If $\a$ and $\b$ are homotopic, then the linking manifolds $L_\a$ and $L_\b$ are cobordant.
\end{prop}


\begin{cor}\label{independenceoflift}
If $\map(P_1,N)$ and $\map(P_2,N)$ are simply connected, then the cobordism class of the linking manifold of $(f_{1,1},f_{2,1})$ only depends on the choice of basepoint $(f_1,f_2)$.
\end{cor}

\begin{pf}
If $\map(X,Y)$ is simply connected, then any two paths in $\map(X,Y)$ sharing the same endpoints are homotopic through paths fixed at their endpoints. The result follows from Proposition \ref{homotopiccobordant}.
$\Box$\end{pf}

\begin{rem}\label{connectivityofmaps}
$\map(X,\RR^k)$ is contractible. In general, it is straightforward to show that $\map(X,Y)$ is $(\conn(Y)-\dim(X))$-connected by considering it as the space of sections of a trivial bundle.
\end{rem}

\begin{rem}
In the recent work of Chernov and Rudyak \cite{cr}, the authors have a way to produce a linking number for a given link from a choice of the unlink in the case where $p_1+p_2+1=n$. It is an element of a cobordism group modulo a group of indeterminacies. Clearly our construction depends also on a choice of path in the space $\map(P_1,N)\times\map(P_2,N)$. In the case where $\map(P_1,N)\times\map(P_2,N)$ is simply-connected, this group of indeterminacies vanishes. This group of indeterminacies also vanishes in other special situations. See especially Section 6 of \cite{cr} for details.
\end{rem}

\begin{rem}\label{independenceofbasepoint}
If we choose another basepoint $(f_1',f_2')\in\link(P_1,P_2;N)$ in the same path component as $(f_1,f_2)$, then the linking manifold remains unchanged. In particular, this is the case if $N$ is connected and $(f_1,f_2)$ and $(f_1',f_2')$ send $P_1$ and $P_2$ to a distinct pair of points in $N$.
\end{rem}

\subsection{Comparison with the classical cases}\label{classicalcase}

Suppose $f_1:S^{p_1}\rightarrow \RR^{p_1+p_2+1}$ and $f_2:S^{p_2}\rightarrow \RR^{p_1+p_2+1}$ are smooth maps with disjoint images. This produces a map $f_1\times f_2:S^{p_1}\times S^{p_2}\rightarrow \RR^{p_1+p_2+1}\times\RR^{p_1+p_2+1}-\Delta_{\RR^{p_1+p_2+1}}$, and by subtraction and retraction, a map $F:S^{p_1}\times S^{p_2}\rightarrow S^{p_1+p_2}$. One classical definition of the \emph{linking number} of $f_1$ and $f_2$ is defined as the degree of this map, which may be computed by computing the cardinality of the inverse image of a regular value $x\in S^{p_1+p_2}$ with signs. The sign convention is as follows. Give an orientation to $S^{p_1}$, $S^{p_2}$, and $S^{p_1+p_2+1}$. Each point $(x_1,x_2)$ in $F^{-1}(x)\subset S^{p_1}\times S^{p_2}$ is assigned the value $+1$ or $-1$ according to whether the derivative $DF:T_{x_1}S^{p_1}\times T_{x_2}S^{p_2}\rightarrow T_xS^{p_1+p_2}$ is orientation preserving or reversing respectively. We will examine this is a slightly more general context where the $p_1$-sphere and the $p_2$-sphere are replaced by smooth closed compact manifolds $P_1$ and $P_2$ of dimensions $p_1$ and $p_2$ respectively. If these manifolds are oriented, then the sign convention is the same as that above, and if either of them is not oriented, then the linking number is an integer mod $2$.\\



Let us first deal with the right way to count the linking number from our perspective. By proposition \ref{countclasses}, $\Om_0^{T\RR^{p_1+p_2+1}-TP_1\times TP_2}(E_{12})$ is $\ZZ$ if $P_1$ and $P_2$ are oriented, and $\ZZ/2$ otherwise.


We have a map $l_2:\hofiber(\link(P_1,P_2;\RR^{p_1+p_2+1})\rightarrow\map(P_1,\RR^{p-1+p_1+1})\times\map(P_2,\RR^{p_1+p_2+1}))\rightarrow\Om C_2(P_1,P_2;\RR^{p_1+p_2+1}))$. Hence for each $$\a\in\hofiber(\link(P_1,P_2;\RR^{p_1+p_2+1})\rightarrow\map(P_1,\RR^{p_1+p_2+1})\times\map(P_2,\RR^{p_1+p_2+1})$$ we get a manifold $L$ of dimension $0$. Taking the induced map on $\pi_0$ and using Proposition \ref{loopcobordism}, we see that $l(\a)$ determines a class in $\Om_0^{TN-TP_1\times TP_2}(E_{12})$.

Next we will show that our definition agrees with the classical one when the basepoint is chosen so that the images of $f_1$ and $f_2$ lie inside disjoint open balls, and that the choice of lift $\a$ is immaterial.\\

Suppose that $\a,\b\in\hofiber(\link(P_1,P_2;\RR^{p_1+p_2+1})\rightarrow\map(S^p,\RR^{p+q+1})\times\map(S^q,\RR^{p+q+1}))$ satisfy $\a(1)=\b(1)$. We write $\a(t)=(f_{1,t},f_{2,t})$ as before, and $\b(t)=(f'_{1,t},f'_{2,t})$. Let $H$ be the straight line homotopy between them. That is, $H(s,t)=(sf_{1,t}+(1-s)f'_{1,t}, sf_{2,t}+(1-s)f'_{2,t})$. Then $H$ gives rise to a cobordism between $L_\a$ and $L_\b$.\\

Consider the composite map

$$
P_1\times P_2\stackrel{f_1\times f_2}{\rightarrow} \RR^{p_1+p_2+1}\times\RR^{p_1+p_2+1}-\Delta_{\RR^{p_1+p_2+1}}\stackrel{d}{\rightarrow} \RR^{p_1+p_2+1}-\{0\}\stackrel{r}{\rightarrow} S^{p_1+p_2}.
$$

As we mentioned, the linking number is the degree of the composed map. As a manifold, it is the inverse image of a regular value $x\in S^{p_1+p_2}$, counted with signs (equivalently, framed) as described above.

Let $L=(r\circ d\circ(f_1\times f_2))^{-1}(x)$. Let $R_x=\{tx\in \RR^{p_1+p_2+1}-\{0\} : t\in(0,\infinity)\}$. $R_x$ is the inverse image of $x$ by $r$, and define $L_1=(d\circ(f_1\times f_2))^{-1}(R_x)$.

\begin{lem}
$L=L_1$ as framed manifolds.
\end{lem}

\begin{pf}
The two sets are equal since $R_x$ is the inverse image of $x$ by $r$. The sign associated to a point in $L$ and the corresponding point in $L_1$ have the same sign provided the orientation of $\RR^{p_1+p_2+1}-\{0\}$ is chosen so that $r$ is an orientation preserving map.
$\Box$\end{pf}

Let $S_d=\{(y,y+tx) | y\in\RR^{p_1+p_2+1},t\in(0,\infinity)\}$. $S_d$ is the inverse image of $R_x$ by the map $d$. It is a $(p_1+p_2+2)$-dimensional set in $\RR^{p_1+p_2+1}\times \RR^{p_1+p_2+1}-\Delta_{\RR^{p_1+p_2+1}}$. Define $L_2=(f\times g)^{-1}(S_d)$.

\begin{lem}
$L_2=L_1$ as framed manifolds.
\end{lem}

\begin{pf}
Since $S_d$ is the inverse image of $R_x$ by $d$, it is clear that $L_1$ and $L_2$ are the same point sets. Again, corresponding points have the same sign provided an orientation of $\RR^{p_1+p_2+1}\times \RR^{p_1+p_2+1}-\Delta_{\RR^{p_1+p_2+1}}$ has been chosen so that $d$ is orientation preserving.
$\Box$\end{pf}

Consider the map $F:P_1\times P_2\times [0,\infinity)\rightarrow N\times N$, where $F(x_1,x_2,t)=(f_1(x_1),f_2(x_2)-tx)$. Define $L_3=F^{-1}(\Delta_N)$.

\begin{lem}
$L_3\cong L_2$ as framed manifolds.
\end{lem}

\begin{pf}
Since the images of $f_1$ and $f_2$ are disjoint, $L_3$ lies away from $P_1\times P_2\times \{0\}$. Hence $L_2$ and $L_3$ are the same point sets counted with signs since they are both the solution sets to the equations $f_1(x_1)=y$ and $f_2(x_2)=y+tx$. The isomorphism between them associates $(x_2,x_2,t)\in L_3$ with $(x_1,x_2)\in L_2$.
$\Box$\end{pf}

Note that by compactness of $P_1$ and $P_2$, there is some $M>0$ such that $L_3\subset P_1\times P_2\times (0,M)$. Hence $f_1(P_1)$ and $f_2(P_2)-tM$ are disjoint and can be separated by enclosing each in a disjoint closed ball. Let these closed disks be denoted $D_{P_1}$ and $D_{P_2}$, with centers $c_{P_1}$ and $c_{P_2}$ respectively. Define $H:P_1\times P_2\times I\rightarrow N\times N$ by the formula

$$H(x_1,x_2,t)= \left\{ \begin{array}{ll}
                        F(x_1,x_2,2Mt) & \mbox{if $0\leq t \leq 1/2$}\\
                         (2-2t)(f_1(x_1),f_2(x_2)-Mx)+(2t-1)(c_{P_1},c_{P_2}) & \mbox{if $1/2\leq t \leq 1$}
                    \end{array}
\right.$$

Note that $H(x_1,x_2,0)=(f_1(x_1),f_1(x_2))$, that $H(x_1,x_2,1)=(c_{P_1},c_{P_2})$, and that $p_1H(x_1,x_2,t)\in\map(P,\RR^{p_1+p_2+1})$ and $p_2H(x_1,x_2,t)\in\map(Q,\RR^{p_1+p_2+1})$ for all $t$. Hence $H(x_1,x_2,t)\in\hofiber(\link(P_1,P_2,\RR^{p_1+p_2+1})\rightarrow\map(P_1,\RR^{p_1+p_2+1})\times\map(P_2,\RR^{p_1+p_2+1}))$, where the homotopy fiber is taken over the constant map $(c_{P_1},c_{P_2})\in\map(P_1,\RR^{p_1+p_2+1})\times\map(P_2,\RR^{p_1+p_2+1})$.

\begin{lem}
$L_3\cong H^{-1}(\Delta_\RR^{p_1+p_2+1})$.
\end{lem}

\begin{pf}
This follows from the observations in the preceding paragraph, together with the observation that $H^{-1}(\Delta_\RR^{p_1+p_2+1})\subset P_1\times P_2\times(0,1/2)$.
$\Box$\end{pf}

\begin{thm}
$L\cong L_3$ as point sets counted with signs, and hence our definition of linking number and the classical one agree.
\end{thm}

\begin{pf}
What remains to check is that $L_3$ is independent of the particular choice of element of $\hofiber(\link(P_1,P_2,\RR^{p_1+p_2+1})\rightarrow\map(P_1,\RR^{p_1+p_2+1})\times\map(P_2,\RR^{p_1+p_2+1}))$. This follows from Proposition \ref{independenceoflift} and Remark \ref{independenceofbasepoint}.
$\Box$\end{pf}

\subsection{Linking in $S^{p_1+p_2+1}$}

By imposing some mild dimensional assumptions we can argue by using our work in the previous section that the classical linking number for linking in spheres agrees with ours. The key observation is that made in Remark \ref{connectivityofmaps}: $\map(X,Y)$ is $(\conn(Y)-\dim(X))$-connected. Since $S^{p_1+p_2+1}$ is $(p_1+p_2)$-connected, $\map(P_1,S^{p_1+p_2+1})$ and $\map(P_2,S^{p_1+p_2+1})$ are simply connected provided $p_1,p_2\geq 1$. In order to use the straight line homotopy as we did before, we need to make sure that every $1$-simplex (homotopy of maps) in $\map(P_1,S^{p_1+p_2+1})$ and $\map(P_2,S^{p_1+p_2+1})$ misses a point (so that they may be considered as maps to $\RR^{p_1+p_2+1}$ and thus we can use the straight line homotopy between them). This will happen provided $p_1+1<p_1+p_2+1$ and $p_2+1<p_1+p_2+1$, or $p_1,p_2\geq 1$.

\section{Cubic Obstructions}\label{linktot2}

In this section we seek to describe the difference between $\link(P_1,\ldots,P_k;N)$ and $\mathcal{T}_2\link(P_1,\ldots,P_k;N)$ in a similar fashion to the way we described the difference between $\link(P_1,\ldots,P_k;N)$ and $\mathcal{T}_1\link(P_1,\ldots,P_k;N)$ in the previous section. That is, we will give a map from $\hofiber(\link(P_1,\ldots,P_k;N)\rightarrow\mathcal{T}_2\link(P_1,\ldots,P_k;N))$ to a cobordism space. As an example, we will use this map to show that the Borromean rings are linked in section \ref{borromeanrings}. We begin by giving a cobordism space model for $\hofiber(\link(P_1,P_2,P_3;N)\rightarrow\mathcal{T}_2\link(P_1,P_2,P_3;N))$. We then go about constructing a manifold in a similar manner to that above.

\subsection{Cobordism model for the cubic stage}\label{S:cubiccobmodel}

Choose a basepoint $({\bf{f}},{\bf{F}},{\bf{H}})\in \mathcal{T}_2\link(P_1,P_2,P_3;N)$. As before ${\bf{f}}=(f_1,f_2,f_3)$, ${\bf{F}}=(F_{12},F_{23},F_{31})$, and ${\bf{H}}=(H_{12},H_{23},H_{31})$. $H_{ij}:P_i\times P_j\times I\rightarrow N\times N$ is a homotopy between $H_{ij}(0)=(f_i,f_j)$ and $H_{ij}(1)=F_{ij}$.

Consider the following diagram $\mathcal{D}$:

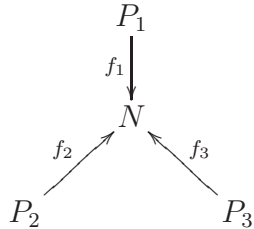
$$
\xymatrix{
&P_1\ar[d]_{f_1}&\\
&N&\\
P_2\ar[ru]^{f_2}&&P_3\ar[lu]_{f_3}\\
}$$

\begin{defn}
Define $E_{123}=\holim(\mathcal{D})$
\end{defn}

Thus a point in $E_{123}$ is a tuple $(x_1,x_2,x_3,\w_1,\w_2,\w_3)$, where $x_i\in P_i$, and $\w_i:I\rightarrow N$ is a path such that $\w(0)=f_i(x_i)$, and $\w_1(1)=\w_2(1)=\w_3(1)$ for $i=1,2,3$. Evidently $E_{123}$ has a map to $P_1\times P_2\times P_3$ given by projection, which we use to pull back $TP_1\times TP_2\times TP_3$ to $E_{123}$. We also have a map $E_{123}\rightarrow N\times N$, defined as follows. There is a map $E_{123}\rightarrow\prod_{i\neq j}\holim(P_i\stackrel{f_i}{\rightarrow
}N\stackrel{f_j}{\leftarrow}P_j)$ given by $(\w_1,\w_2,\w_3)\mapsto(\w_1\cdot\w_2^{-1}, \w_2\cdot\w_3^{-1}, \w_3\cdot\w_1^{-1})$, where $\cdot$ indicates path multiplication. If we reparametrize the resulting paths to be defined on $[-1,1]$ then the map to $E_{123}\rightarrow N\times N$ cam be defined as the composition of this map with the evaluation $(\w_1\cdot\w_2^{-1}, \w_2\cdot\w_3^{-1}, \w_3\cdot\w_1^{-1})\mapsto(\w_1\cdot\w_2^{-1}(0),\w_2\cdot\w_3^{-1}(0))$, and we use this map to pull back $TN\times TN$ to $E_{123}$. We use this to define our cobordism space.

\begin{defn}
Define $C_3(P_1,P_2,P_3;N)=C_\bullet^{TN\times TN-TP_1\times TP_2\times TP_3}(E_{123})$.
\end{defn}

In this case, $d=p_1+p_2+p_3-2n$, but for indexing purposes, it is more convenient to shift this down by $2$ (this could be done honestly by subtracting twice the trivial $1$-dimensional bundle, but that is rather cumbersome notationally). Thus, a $0$-simplex in this space will be represented by a manifold of dimension $p_1+p_2+p_3-2n+2$. Note that $C_3$ is equivalent to $QT(E_{123};TN\times TN-TP_1\times TP_2\times TP_3)$ by the Pontryagin-Thom construction.

\subsection{Whitney Circles and Whitney Disks Terminology}


Proposition \ref{t2model} tells us that we may think of an element of $\mathcal{T}_2\link(P_1,P_2;N)$ as a triple $((f_1,f_2), F_{12}, H_{12})$ where $f_i:P_i\rightarrow N$, $F_{12}:P_1\times _2\rightarrow N\times N-\Delta_N$, and $H_{12}:P_1\times P_2\times I\rightarrow N\times N$ is a homotopy from $f_1$ to $f_2$. Suppose $H_{12}$ is transverse to $\Delta_N$. Then $D_{12}=H_{12}^{-1}(\Delta_N)$ is a $(p_1+p_2+1-n)$-dimensional submanifold of $P_1\times P_2\times I$ which we call the \emph{Whitney circle}. The reason for this terminology is that in case $D_{12}$ is $1$-dimensional, the images under $f_i$ of the projections of $D_{12}$ to $P_i$ for $i=1,2$ form a bent circle in $N$. This circle bounds a disk in $N$ as follows. For $(x_1,x_2,t_0)\in D_{12}$ we give a path in $N$ from $f_1(x_1)$ to $f_2(x_2)$ by going from $f_1(x_1)=p_1H_{12}(x_1,x_2,0)$ to $p_1H_{12}(x_1,x_2,t_0)=p_2H_{12}(x_1,x_2,t_0)$ by letting $t$ run from $0$ to $t_0$, and then to $f_2(x_2)$ by letting $t$ run from $t_0$ back down to $0$. We call this disk (and its higher-dimensional counterpart) the \emph{Whitney disk}. Our obstruction manifold measures the intersections of the Whitney circles, and the intersections of the manifolds themselves with the Whitney disks. Compare an identical description in \cite{st} in the case of $2$-spheres in a $4$-manifold.

\subsection{Construction of the obstruction manifold}

For the reader's sake we will construct the obstruction manifold in the case $k=3$, so the space under consideration is $\link(P_1,P_2,P_3;N)$. For higher $k$, one constructs $k\choose 3$ manifolds in the manner described below. There are two types of intersections which form the obstruction manifold. The first is the intersections of $P_j$ with the generalized Whitney disks for the pair $(P_h,P_i)$, where $h,i,j$ are distinct. The second is the intersections of the generalized Whitney circles of the pair $(P_h,P_i)$ with that for $(P_i,P_j)$, again for $h,i,j$ distinct.


\subsubsection{Intersection of the Whitney disk with the $P_i$}

Let $T_1=\{(s,t) | 0\leq s\leq t \leq 1\}$. For $(h,i,j)\in S$ consider the maps $\Phi_{1,(hi)j}:P_1\times P_2\times P_3\times T\rightarrow N^4$ given by

\begin{equation}
\Phi_{1,(hi)j}(x_1,x_2,x_3,s,t)=(f_j(x_j), p_1H_{hi}(x_h,x_i,s), H_{hi}(x_h,x_i,t)) 
\end{equation}

Let $T_2=\{(s,t) | 0\leq t\leq s \leq 1\}$. For $(h,i,j)\in S$ consider the maps $\Phi_{2,(hi)j}:P_1\times P_2\times P_3\times T\rightarrow N^4$ given by

\begin{equation}
\Phi_{2,(hi)j}(x_1,x_2,x_3,s,t)=(H_{hi}(x_h,x_i,s), p_2H_{hi}(x_h,x_i,t), f_j(x_j)) 
\end{equation}


\begin{defn}
For $a=1,2$, define $X_{a,(hi)j}=\Phi_{a,(hi)j}^{-1}(\Delta_N\times\Delta_N)$.
\end{defn}

\begin{prop}\label{Xmanifold}
For $a=1,2$, if $\Phi_{a,(hi)j}$ and its restrictions to $\del_1 (P_1\times P_2\times P_3 \times T_a)$ and $\del_2 (P_1\times P_2\times P_3 \times T_a)$ are transverse to $\Delta_N\times\Delta_N$, then $X_{a,(hi)j}$ is a compact $(p_1+p_2+p_3+2-2n)$-dimensional manifold with at most $2$-strata, and stable normal bundle $TN\times TN- TP_1\times TP_2\times TP_3$.
\end{prop}

\begin{pf}
Transversality gives an isomorphism $TX_{a,(hi)j}\oplus TN\oplus TN\rightarrow TP_1\oplus TP_2\oplus TP_3\oplus \e^2$.
$\Box$\end{pf}

The boundary of $X_{a,(hi)j}$ naturally decomposes in to three parts according to whether $t=0, s=0$, or $t=s$. $\del X_{a,(hi)j}=\del_tX_{a,(hi)j}\cup \del_sX_{a,(hi)j}\cup \del_{t=s}X_{a,(hi)j}$. Also, $\del_{t=s}X_{1,(hi)j}=\del_{t=s}X_{2,(hi)j}$. The $2$-stratum of $X_{a,(hi)j}$ consists of the points $(p_1,p_2,p_3,0,0)$ such that $f_1(p_1)=f_2(p_2)=f_3(p_3)$, a compact $(p_1+p_2+p_3-2n)$-dimensional manifold we will call $T$.

\subsubsection{Intersections of the Whitney circles}

For $(h,i,j)=(1,2,3), (2,3,1), (3,1,2)$, consider the maps $\Phi_{hij}:P_1\times P_2\times P_3\times I\times I\rightarrow N^4$ given by

\begin{equation}
\Phi_{hij}(p_1,p_2,p_3,s,t)=(H_{hi}(p_h,p_i,s),H_{ij}(p_i,p_j,t))
\end{equation}


\begin{defn}
$W_{hij}=\Phi_{hij}^{-1}(\Delta_N\times \Delta_N)$.
\end{defn}

\begin{prop}\label{Wmanifold}
If $\Phi_{hij}$ and its restrictions to $\del_1(P_1\times P_2\times P_3 \times I\times I)$ and $\del_2 (P_1\times P_2\times P_3 \times I\times I)$ are transverse to $\Delta_N\times\Delta_N$, then $W_{hij}$ is a compact $(p_1+p_2+p_3+2-2n)$-dimensional manifold with boundary, and stable normal bundle $TN\times TN-TP_1\times TP_2\times TP_3$.
\end{prop}

\begin{pf}
Transversality gives an isomorphism $TW_{hij}\oplus TN\oplus TN\rightarrow TP_1\oplus TP_2\oplus TP_3\oplus \e^2$.
$\Box$\end{pf}

The boundary $\del W_{hij}$ occurs when either $s=0$ or $t=0$, and gives a decomposition $\del W_{hij}=\del_s W_{hij}\cup \del_t W_{hij}$. The $2$-stratum of $W_{hij}$ is the set $T$ described above.

\subsubsection{Forming the obstruction manifold}

The following lemma tells us how the boundaries of the $W_{hij}$ and $X_{a,(hi)j}$ fit together, and can be verified by noting that the equations which define these manifolds are the same when the appropriate value of $t$ or $s$ is $0$.

\begin{lem}
For $(h,i,j)\in S$, there are diffeomorphisms of manifolds with boundary $\del_tW_{hij}\cong \del_t X_{2,(hi)j}$, $\del_sW_{hij}\cong \del_sX_{1,(ij)h}$, and $\del_{s=t}X_{1,(hi)j}\cong \del_{s=t}X_{2,(hi)j}$.
\end{lem}

We define $\del X_{(hi)j}=\del_{s=t}X_{1,(hi)j}$, and arrange this into the following diagram:

$$
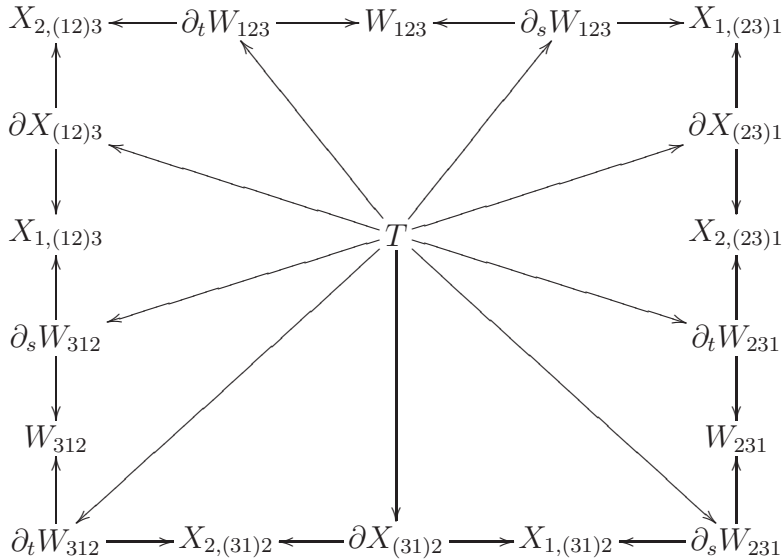
\xymatrix{
X_{2,(12)3} & \del_tW_{123}\ar[l]\ar[r] & W_{123} & \del_sW_{123}\ar[l]\ar[r] & X_{1,(23)1} \\
 \del X_{(12)3}\ar[u]\ar[d] &  &  &  & \del X_{(23)1}\ar[u]\ar[d]  \\
X_{1,(12)3} &  & T\ar[luu]\ar[ruu]\ar[llu]\ar[rru]\ar[lld]\ar[rrd]\ar[llddd]\ar[rrddd]\ar[ddd] &  & X_{2,(23)1} \\
 \del_sW_{312}\ar[u]\ar[d] &  &  &  & \del_tW_{231}\ar[u]\ar[d] &  \\
W_{312} &  &  &  & W_{231} &  &  \\
\del_tW_{312}\ar[u]\ar[r] & X_{2,(31)2} & \del X_{(31)2}\ar[l]\ar[r] & X_{1,(31)2} &  \del_sW_{231}\ar[l]\ar[u] &  \\
}$$

Let $\mathcal{I}$ denote the underlying category for this diagram, and let $\mathcal{Z}:\mathcal{I}\rightarrow\mbox{ Spaces}$ be the functor which gives the diagram above. 

\begin{defn}
Define $Z=\colim_{\mathcal{I}}\mathcal{Z}$, and let $Z'=\hocolim_{\mathcal{I}}\mathcal{Z}$.
\end{defn}

Note that the canonical map $a:Z'\rightarrow Z$ is a homotopy equivalence. We wish to show how to piece together the bundle isomorphisms given by transversality in Propositions \ref{Xmanifold} and \ref{Wmanifold} to make a bundle isomorphism over $Z'$. In order to do this, we need to compare the isomorphisms along the common boundary they define and invoke Propositions \ref{gluebundles1} and \ref{gluebundles2}. In order to apply these propositions, we need to first describe the map $Z\rightarrow E_{123}$.

\subsection{The map $z:Z\rightarrow E_{123}$}

The map $z:Z\rightarrow E_{123}$ comes in several pieces, because $Z$ itself is defined in several pieces. Note that, by composition, this induces a map $z':Z'\rightarrow E_{123}$. A point in $Z$ is a tuple $(x_1,x_2,x_3,s,t)$ which is in one of the $W_{hij}$ or $X_{a,(hi)j}$. If $(x_1,x_2,x_3,s,t)\in W_{hij}$, then the corresponding point in $E_{123}$ has

$$\w_h(r_h)= \left\{ \begin{array}{ll}
                        p_1H_{hi}(x_h,x_i,2r_hs) & \mbox{if $0\leq r_h \leq 1/2$}\\
                        p_2H_{hi}(x_h,x_i,2(1-r_h)s) & \mbox{if $1/2\leq r_h \leq 1$}
                    \end{array}
\right.,$$

$\w_i(r_i)=f_i(x_i)$, and

$$\w_j(r_j)= \left\{ \begin{array}{ll}
                        p_2H_{ij}(x_i,x_j,2r_jt) & \mbox{if $0\leq r_j \leq 1/2$}\\
                        p_1H_{ij}(x_i,x_j,2(1-r_j)t) & \mbox{if $1/2\leq r_j \leq 1$}
                    \end{array}
\right.,$$

for all $r_i\in [0,1]$. If $(x_1,x_2,x_3,s,t)\in X_{1,(hi)j}$, then

$\w_h(r_h)=p_1H_{hi}(x_h,x_i,r_hs)$

$$\w_i(r_i)= \left\{ \begin{array}{ll}
                        p_2H_{hi}(x_h,x_i,2r_it) & \mbox{if $0\leq r_i \leq 1/2$}\\
                        p_1H_{hi}(x_h,x_i,2(1-r_i)t+2(r_i-1/2)s) & \mbox{if $1/2\leq r_j \leq 1$}
                    \end{array}
\right.,$$ and

$w_j(r_j)=f_j(x_j)$ for all $r_j\in [0,1]$. Finally, if $(x_1,x_2,x_3,s,t)\in X_{2,(hi)j}$, then

$$\w_h(r_h)= \left\{ \begin{array}{ll}
                        p_1H_{hi}(x_h,x_i,2r_hs) & \mbox{if $0\leq r_j \leq 1/2$}\\
                        p_2H_{hi}(x_h,x_i,2(1-r_h)s+2(r_h-1/2)t) & \mbox{if $1/2\leq r_j \leq 1$}
                    \end{array}
\right.,$$

$\w_i(r_i)=p_2H_{hi}(x_h,x_i,r_it)$, and $\w_j(r_j)=f_j(x_i)$ for all $r_j\in [0,1]$. By inspection, these maps agree on the intersections of the $W_{hij}$ with the $X_{a,(hi)j}$. Hence we have a map $Z\rightarrow E_{123}$. Now we are going to produce a path in $C_{3}(P_1,P_2,P_3;N)$ from an element $\a$ of $$\hofiber(\link(P_1,P_2,P_3;N)\rightarrow\mathcal{T}_2\link(P_1,P_2,P_3;N))\rightarrow \Om C_{3}(P_1,P_2,P_3;N).$$

Write $\a(u)=(H'_{12}(u),H'_{23}(u),H'_{31}(u))$, where $u\in I$, $\a(0)=(H_{12},H_{23},H_{31})\in\mathcal{T}_2\link(P_1,P_2,P_3;N)$, and $\a(1)=(H'_{12},H'_{23},H'_{31})$. Here $H'_{ij}:P_i\times P_j\times I\rightarrow N\times N$ satisfies $H'_{ij}=f'_i\times f'_j$ for some $(f'_1,f'_2,f'_3)\in\link(P_1,P_2,P_3;N)$. Since we can think of $\a$ as a parametrized family of maps in $\mathcal{T}_2\link(P_1,P_2,P_3;N)$, we can apply the construction of $Z$ to this family. What remains is to give the map to $E_{123}$, which means associating three paths with the point. For the sake of brevity, let us do this in the specific case where $(x_1,x_1,x_3,s,t,u)\in W_{123}$. This means that $H_{12}(x_1,x_2,s,u)\in\Delta_N$ and $H_{23}(x_2,x_3,t,u)\in\Delta_N$. Then, for instance, we define

$$\w_1(r_1)= \left\{ \begin{array}{ll}
                        p_1H_{12}(x_1,x_2,2r_1s,2r_1u) & \mbox{if $0\leq r_j \leq 1/2$}\\
                        p_2H_{12}(x_1,x_2,2(1-r_1)s,2(1-r_1)u) & \mbox{if $1/2\leq r_j \leq 1$}
                    \end{array}
\right..$$

The others are defined similarly.

\subsection{Gluing the bundle data}

For vector spaces $V,W$, let $L_{\max}(V,W)$ denote the space of linear maps $V\rightarrow W$ of maximal rank. There is an isomorphism $L_{\max}(V,W)\rightarrow L_{\max}(W,V)$ given by sending a linear transformation to its transpose. Note that when $\dim(V)\leq \dim(W)$, space $L_{\max}(V,W)$ is the Stiefel manifold, which is well-known to be $(\dim(W)-\dim(V)-1)$-connected. We wish to consider bundles of spaces $L_{\max}(V,W)$, as introduced in Definition \ref{maxrankbundlemaps}. 

\begin{lem}\label{derivativeshomotopic}
The derivatives $D\Phi_{hij}$ and $D\Phi_{a,(hi)j)}$ are homotopic where it makes sense to compare them.
\end{lem}

\begin{pf}
Consider the space $L_{\max}(TP\times\e^2,TN\times TN)$ as a space over $\del_tW_{hij}\cong \del_tX_{2,(hi)j}$. The maps $D\Phi_{hij}$ and $D\Phi_{2,(hi)j)}$ define sections of these bundles, and we wish to show they are homotopic.

Take an open cover $\mathcal{U}$ of $\del_tW_{hij}$ consisting of contractible open sets $U$ so that the bundle $L_{\max}(TP\times\e^2,TN\times TN)$ is trivial over each $U\in\mathcal{U}$. Now $L_{\max}(TP\times\e^2,TN\times TN)$ is $(p_1+p_2+p_3-2n+1)$-connected, and if  $p_1+p_2+p_3+2-2n\leq 0$, there is nothing to prove. Otherwise, since $(p_1+p_2+p_3-2n+1)\geq 0$, two such maps are always homotopic. It follows that $D\Phi_{hij}$ and $D\Phi_{2,(hi)j}$ are homotopic, by Remark \ref{sectionslinear}. The same proof shows that the restrictions of $D\Phi_{hij}$ and $D\Phi_{1,(ij)h}$ and $D\Phi_{1,(ij)h}$ and $D\Phi_{2,(ij)h}$ are homotopic where it makes sense to compare them. As for the $2$-strata $T$, we must compare these chosen homotopies on $T$ and ensure that the composed homotopy is homotopic to a constant homotopy. If $p_1+p_2+p_3+2-2n\leq 1$ there is nothing to prove, since the $2$-strata is empty. Otherwise $(p_1+p_2+p_3-2n+1)\geq 1$ and the space $L_{\max}(TP\times\e^2,TN\times TN)$ is $\geq1$-connected.
$\Box$\end{pf}

\begin{lem}\label{piecetogetherbundledata}
There is an isomorphism $TZ\oplus TN\oplus TN\rightarrow TP_1\oplus TP_2\oplus TP_3\oplus \e^2$.
\end{lem}


\begin{pf}
Let $P=P_1\times P_2\times P_3$, and identify $\nu(\Delta_N\times \Delta_N\subset N^4)$ with $TN\times TN$. 
For $(h,i,j)\in S$, the derivatives $D\Phi_{hij}:TP\times \e^2\rightarrow TN^4$ give sections of the bundles $L_{\max}(TP\times\e^2,TN\times TN)$ over $W_{hij}$. The same is true of $D\Phi_{a,(hi)j}:TP\times\e^2\rightarrow TN\times TN$ and the manifolds $X_{a,(hi)j}$. These maps are homotopic by Lemma \ref{derivativeshomotopic}. By Propositions \ref{gluebundles1} and \ref{gluebundles2}, they combine to give an element of $\hocolim_{\mathcal{I}}L_{\max}(z'^\ast TP\times\e^2,z'^\ast TN\times TN)$, which is a space over $Z'$. Pulling back by a homotopy inverse to $a:Z'\rightarrow Z$, we obtain an element $\Phi\in L_{\max}(z^\ast TP\times \e^2,z^\ast TN\times TN)$. The fiberwise kernel of $\Phi$ is, by construction, $TZ$ (see Lemma \ref{atmost2stratasmoothmanifold}), and this gives the desired isomorphism.
$\Box$\end{pf}

\begin{thm}
$Z$ can be given the structure of a smooth closed compact manifold of dimension $p_1+p_2+p_3+2-2n$, and there is an isomorphism $TZ\oplus TN\oplus TN\cong TP_1\oplus TP_2\oplus TP_3\oplus \e^2$.
\end{thm}

\begin{pf}
That $Z$ is smooth follows from Proposition \ref{atmost2stratasmoothmanifold}, since the diffeomorphisms between the $1$-strata of the pieces which form $Z$ all restrict to the identity on $T$. Note that the canonical map $a:Z'\rightarrow Z$ is a homotopy equivalence. By Lemma \ref{piecetogetherbundledata}, we pull back the bundle isomorphism over $Z'$ by a homotopy inverse to $z$ to give the desired bundle isomorphism.
$\Box$\end{pf}

What we have shown is that there is a map from the subcomplex of $0$-simplices of $$\hofiber(\link(P_1,P_2,P_3;N)\rightarrow\mathcal{T}_2\link(P_1,P_2,P_3;N))$$ such that the above transversality conditions are met, to the $0$-simplices of $\Om C_3(P_1,P_2,P_3;N)$. Getting a map of simplicial sets is easy, since the above works for families of maps too, and since the subcomplex of transverse maps is homotopy equivalent to the full complex. Thus we have proven Theorem \ref{t2mapofspaces}, restated below.


\begin{thm}
For a basepoint in the image of $\link(P_1,P_2,P_3;N)$, there is a map of spaces $$l_3:\hofiber(\link(P_1,P_2,P_3;N)\rightarrow\mathcal{T}_2\link(P_1,P_2,P_3;N))\rightarrow \Om C_{3}(P_1,P_2,P_3;N).$$
\end{thm}

In the event the basepoint is a point in $\mathcal{T}_2\link(P_1,P_2,P_3;N)$, we construct the manifold $Z$ from this point as above, and the relevant cobordism space is $\Om_Z C_{3}(P_1,P_2,P_3;N)$, the space of nullcobordisms of $Z$. Compare Theorem \ref{t1mapofspaces}.

We conjecture the following analog of Theorem \ref{c2connectivity}, to which we will dedicate a future paper:

\begin{conj}\label{c3connectivity}
$l_3$ is $(3n-\max\{2p_1+p_2+p_3, p_1+2p_2+p_3, p_1+p_2+2p_3\}-5)$-connected.
\end{conj}




\subsection{Example: The Borromean Rings}\label{borromeanrings}


As an example, we show that the Borromean rings are not homotopic to the unlink. We proceed by finding a path from image of the Borromean rings to the unlink in $\mathcal{T}_1\link(S^1,S^1,S^1;\RR^3)$, lifting this path to a path in $\mathcal{T}_2\link(S^1,S^1,S^1;\RR^3)$, and then constructing the obstruction manifold for this family, which turns out to be a single point. 

Let $B=(f_1,f_2,f_3)\in\link(S^1,S^1,S^1;\RR^3)$ be the Borromean rings, where $f_1(x)=(0,2\cos x,\sin x)$, $f_2(y)=(\cos y,0,2\sin y)$, and $f_3(z)=(2\cos z, \sin z, 0)$. Let the unlink be represented by the triple $U=(f_1,e_2,f_3)\in\link(S^1,S^1,S^1;\RR^3)$, where $e_2(y)=(\cos y, 4, \sin y)$. The images of  $B$ and $U$ in $\mathcal{T}_1\link(S^1,S^1,S^1;\RR^3)$ are clearly homotopic by the path $L(t)=(f_1,H_1(y,t),f_3)$, where $H_1(y,t)=(\cos y, t, 2\sin y)$ for $t\in[0,4]$. We lift this path to the path in $\mathcal{T}_2\link(S^1,S^1,S^1;\RR^3)$ given by the triple $(f_1\times f_3, H_1\times f_1, f_3\times H_3)$, where $H_3(y,t)=(\cos y, t, 2\sin y -(t-2)^2+4)$ for $t\in[0,4]$ (we have altered the order in which we represent tuples in $\mathcal{T}_2$ that we used in Section \ref{S:cubiccobmodel}, and the $\mathcal{T}_2$ we mean here is the image of an inner automorphism of $\mathcal{T}_2$ induced by a permutation of $\{1,2,3\}$).


In order to construct the obstruction from this family, we must make the double point construction, and form the Whitney disks. The double points of the family $L(t)$ come from the intersection of $f_3$ with $H_1$, which is clearly a pair of points $(y,z,t)=(\pi/2,\pi/3,\sqrt{3}/2),(3\pi/2,2\pi/3,\sqrt{3}/2)$, and this intersection is transverse. To make a nullcobordism of this manifold, we define a function $G:S^1\times S^1\times [0,4]\times [0,1]\rightarrow \RR^3\times\RR^3$ given by $G(y,z,t,s)=(f_3(z), sH_3(y,t)+(1-s)H_1(y,t))$.

\begin{prop}
$G$ is transverse to $\Delta_{\RR^3}$ and $G^{-1}(\Delta_{\RR^3})\subset S^1\times S^1\times[0,4]\times [0,1]$ is an arc which embeds in each copy of $S^1$ under the projections to the $S^1$ factors. 
\end{prop}

\begin{pf}
Let $H(y,t,s)=sH_3(y,t)+(1-s)H_1(y,t)$. To see that $G$ is transverse to the diagonal, it is equivalent to check that the images of $f_3$ and $H$ intersect transversely in $\RR^3$. The derivative

$$DH=\left( \begin{array}{ccc}
-\sin y & 0 & 0 \\
0 & 1 & 0 \\
2\cos y & -2s(t-2) & -(t-2)^2+4 \\
\end{array}
\right)$$

has determinant $(-(t-2)^2+4)\sin y$, which vanishes if and only if either $t=0,4$ or $y=0,\pi$. Since there are no solutions to $f_3(z)=H(y,0,s)$, we only have to deal with the possibility of $y=0,\pi$. In this case, $z=\pi/3, 2\pi/3$, and by inspection we find that for these values of $z$, $Df_3$ and the last two columns of $DH$ span $\RR^3$.\\

From transversality it follows that $D=G^{-1}(\Delta_{\RR^3})$ is a $1$-dimensional manifold with boundary. It is the solutions to the following equations:

\begin{eqnarray*}
2\cos z = \cos y\\
\sin z = t\\
0 = 2\sin y+s(-(t-2)^2+4)
\end{eqnarray*}

The first two equations imply that $z\in[\pi/3,2\pi/3]$ and $t\in[\sqrt{3}/2,1]$. Moreover, for each such $(z,t)$ there is a unique $(y,s)$ which solves all three, since the third equation implies that $2\sin y \leq 0$, and we can solve for $s$. Let $p_i:S^1\times S^1\times[0,4]\times[0,1]\rightarrow S^1$ for $i=1,2$ be the projections. We have shown that $p_2(D)=[\pi/3,2\pi/3]$ and that the restriction of $p_2$ to $D$ is an embedding. It is also clear that $p_1(D)=[\pi,2\pi]$ and that the projection $p$ of $D$ onto $S^1\times I$, given by $p(y,z,t)=(y,t)$, is an embedding (in fact, the restriction of $p_1$ itself is an embedding, though we will not need this).
$\Box$\end{pf}

The Whitney circle $C=f_3(p_2D)\cup H_1(pD)$ is a bent circle in $\RR^3$. Let $\pi_i:\RR^3\times\RR^3\rightarrow\RR^3$ be the projections for $i=1,2$. $C$ bounds a disk $W$ given by $G$ as follows: For each $(y,z,t,s)\in D$, $G$ gives rise to a path in $\RR^3$ from $f_3(z)=\pi_1G(y,z,t,0)$ to $H_1(y,t)=\pi_2G(y,z,t,0)$ because $\pi_1G(y,z,t,s)=\pi_2G(y,z,t,s)$. The path is defined by

$$\gamma(\s)= \left\{ \begin{array}{ll}
\pi_1G(y,z,t,2\s s) & \mbox{if $\s\in[0,1/2]$}\\
\pi_2G(y,z,t,(2-2\s)s) & \mbox{if $\s\in[1/2,1]$}
\end{array}\right.
$$

The union of these paths forms the disk $W$ bounding $C$ in $\RR^3$, and it is clearly disjoint from $f_3(S^1)$ and $H_1(S^1)\times [0,4]$ except along its boundary. It is clear that $f_2(S^1)$ intersects $W$ transversely in a single point, since $W$ is the union of straight lines between $z$ and $(y,t)$ for each $(y,z,t,s)\in D$ (and for each $y$ there is a unique $z$ and vice-versa, since the projections of $D$ to the $S^1$ factors are embeddings). It follows that $B$ cannot be link homotopic to $U$, for a link homotopy between them would clearly produce an empty manifold, which is not cobordant to a point.




\section{Acknowledgments}

The author would like to thank Greg Arone, Ryan Budney, Vladimir Chernov, Ralph Cohen, Tom Goodwillie, and Steve Kerckhoff for helpful conversations. He would also like to thank Don Stanley, who first posed the question to the author about whether the linking number could be seen using calculus of functors. This research was partly supported by NSF grant DMS-0402822.




\end{document}